\documentclass{article}

\usepackage[margin=1in]{geometry}
\usepackage{amsmath,amsthm,amssymb,amsfonts}
\usepackage{graphicx}
\usepackage{setspace}
\usepackage{enumitem}   
\graphicspath{{pictures/}}
\usepackage{amsmath,amssymb}
\usepackage{hyperref}

\newtheorem{thm}{Theorem}
\newtheorem{prop}{Proposition}
\newtheorem{mydef}{Definition}

\newtheorem{remark}{Remark}

\title{Geography of Genus $2$ Lefschetz Fibrations}
\author{Kai Nakamura}
\date{\today}

\begin{document}
\maketitle

\begin{abstract}
Questions of geography of various classes of $4$-manifolds have been a central motivating question in $4$-manifold topology.
Baykur and Korkmaz asked which small, simply connected, minimal $4$-manifolds admit a genus $2$ Lefschetz fibration.
They were able to classify all the possible homeomorphism types and realize all but one with the exception of a genus $2$ Lefschetz fibration on a symplectic $4$-manifold homeomorphic, but not diffeomorphic to $3 \mathbb{CP}^2 \# 11\overline{\mathbb{CP}}^2$.
We give a positive factorization of type $(10,10)$ that corresponds to such a genus $2$ Lefschetz fibration.
Furthermore, we observe two restrictions on the geography of genus $2$ Lefschetz fibrations, we find that they satisfy the Noether inequality and a BMY like inequality.
We then find positive factorizations that describe genus $2$ Lefschetz fibrations on simply connected, minimal symplectic $4$-manifolds for many of these points.
\end{abstract}
\section{Introduction}
One of the primary objectives of symplectic $4$-manifold topology is to construct simply connected symplectic $4$-manifolds.
A fruitful approach is to use Donaldson and Gompf's characterization of symplectic manifolds as those that admit a Lefschetz pencil \cite{donaldson1999lefschetz,gompf2004toward}.
A genus $g$ Lefschetz fibration can be described combinatorially as a positive factorization, a sequence of Dehn twists that compose to the identity in the mapping class group of a genus $g$ surface.
This idea was used to great effect by Baykur and Korkmaz to construct many positive factorizations that describe Lefschetz fibrations on minimal, simply connected symplectic $4$-manifolds that are homeomorphic, but not diffeomorphic to blow ups of $\mathbb{CP}^2$ and $3\mathbb{CP}^2$ \cite{baykur2017small}.
Furthermore, they were able to classify the possible homeomorphism types for simply connected, minimal small genus $2$ Lefschetz fibrations and this construction realizes all but one which would correspond to a positive factorization of type $(10,10)$.
We construct such a positive factorization by the use of lantern substitutions.
We can then control the minimality of the corresponding Lefschetz fibration by taking advantage of a theorem of Endos and Gurtas, that a lantern substitution in a positive factorization corresponds to a rational blow down of the Lefschetz fibrations \cite{endo2010lantern}.
This corresponds to the $(x,y) = (2,8)$ case of the following theorem which is our main result
\begin{thm}
If $(X,f)$ is a genus $2$ Lefschetz fibration, then the following inequality holds
\[2\chi_h(X)-6 \le c_1^2(X) \le 6\chi_h(X)-3\]
The latter inequality can be made strict if $X$ is simply connected.
For $(x,y)$ such that $2x-6 \le y \le 5.5x-3$ with $y \ge 0$, there exists simply connected, minimal genus $2$ Lefschetz fibrations such that $\chi_h = x$ and $c_1^2 = y$.
\end{thm}

These Lefschetz fibrations will all be given by explicit positive factorizations. 
This theorem represents a strong generalization of the work of Baykur and Korkmaz.
Combined with the results of Sato on the geography of non-minimal genus $2$ Lefschetz fibrations, this gives a fairly complete picture of the geography of genus $2$ Lefschetz fibrations \cite{sato20082}.
The restrictions on the geography will follow from relating the type $(n,s)$ of a genus $2$ Lefschetz fibration $(X,f)$ to the topological invariants of $X$.
Then we will show how Baykur and Korkmaz's work will generalize to the case of $2x-6 \le y \le 5x - 3$ with $y \ge 0$.
Then we will use our construction of a $(10,10)$ positive factorization to construct the remainder.

Part of these results originally appeared in the author's undergraduate thesis  \cite{nakamurathesis}.
The construction of a $(10,10)$ positive factorization, restrictions on the geography, and realization of simply connected, minimal genus $2$ Lefschetz fibrations below the $BK$-line, and a positive factorization of type $(24,38)$ were all originally constructed in the thesis.
\section{Background and Prior Results}
\subsection{Lefschetz Fibrations and Positive Factorizations}
\begin{figure}
\centering
\includegraphics[scale = 0.77]{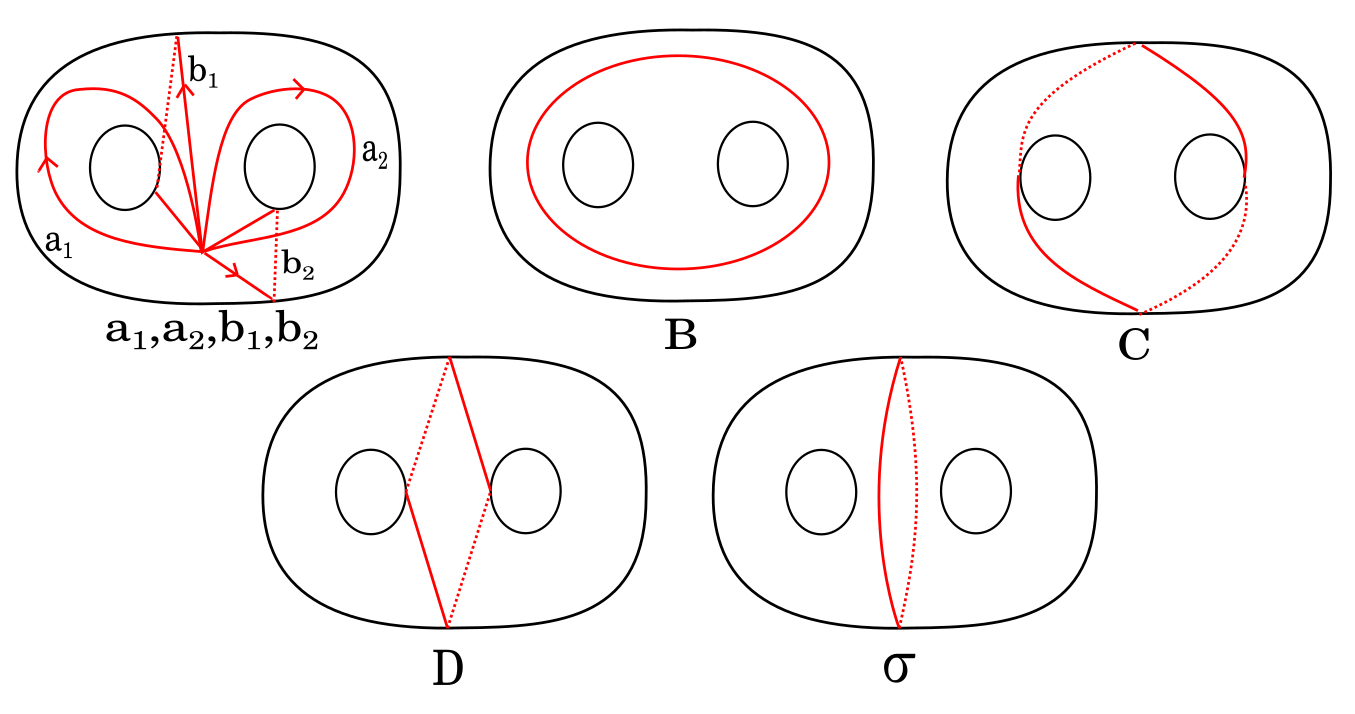}
\caption[matsumoto and generators]{The Generators of $\pi_1(\Sigma_2)$ and the curves $B,C,D,\sigma$ coming from Matsumoto's fibration}
\label{fig:matsumoto}
\end{figure}
Donaldson and Gompf provided a topological characterization of symplectic manifolds as those manifolds that admit a Lefschetz pencil \cite{donaldson1999lefschetz,gompf2004toward}.
The base points of a Lefschetz pencil can be blown up to get a Lefschetz fibration
\begin{mydef}
A Lefschetz fibration $(X,f)$ on a closed, smooth, oriented $4$-manifold $X$ is a smooth surjective map $f: X \rightarrow S^2$ whose critical locus consists of finitely many points $p_i$ and there exist local coordinates at each $p_i$ and $f(p_i)$ (which agree with the orientations on $X$ and $S^2$) where the map takes the form $f(z_1,z_2) = z_1z_2$
\end{mydef}
The non-singular fibers of a Lefschetz fibration will be a surface of some genus $g$ and we call these genus $g$ Lefschetz fibration.
If the complement of a singular fiber is connected, we call the fiber irreducible and otherwise we call the singular fiber reducible.
We can take the critical points to lay on disjoint fibers.
Then each singular fiber will take the form of the non-singular fiber with a simple closed curve $c_i$ collapsed to the critical point $p_i$.
The monodromy of the Lefschetz fibration around one of these fibers will be the Dehn twist $T_{c_i}$ and the composition of all of these Dehn twists will be the identity.
This is called a positive factorization of the identity and we get one for every Lefschetz fibration.

The genus $g$ mapping class group $\Gamma_g$ is the group of homeomorphisms of the genus $g$ surface $\Sigma_g$ up to isotopy.
There is a correspondence between genus $g$ Lefschetz fibrations and positive factorizations of the identity in to Dehn twists in $\Gamma_g$.
This means that if we are given some positive factorization 
\[T_{c_n} \dots T_{c_1} = Id\]
Then there is some genus $g$ Lefschetz fibration with this positive factorization as its monodromy.
This correspondence is up to equivalence of Hurwitz moves and global conjugations.
A Hurwitz move is replacing some $T_{c_{i+1}}T_{c_i}$ in a positive factorization with $T_{T_{c_{i+1}}(c_i)}T_{c_{i+1}}$, these will equal as mapping classes since $T_{T_{c_{i+1}}(c_i)} = T_{c_{i+1}}T_{c_i}T_{c_{i+1}}^{-1}$.
Then global conjugation is when we take some mapping class $\phi$ and replace each $c_i$ with $\phi(c_i)$.
If $W$ is a positive factorization, then we denote by $W^\phi$ the positive factorization coming from global conjugation of $W$ by $\phi$.
When the genus two or larger, then we have a true one-to-one equivalence between positive factorizations of the identity in the genus $g$ mapping class group up to equivalence and genus $g$ Lefschetz fibrations up to isomorphism.

If we have a positive factorization for a genus $g$ Lefschetz fibration $(X,f)$, we can read off many of the topological invariants from this positive factorization.
If the Lefschetz fibration has positive factorization $T_{c_n} \dots T_{c_1}$, then we have a presentation for the fundamental group of $X$ as the quotient of the fiber by the fundamental group of the non-singular fiber by the subgroup normally generated by $c_1,\dots, c_n$
\[\pi_1(X) \cong \pi_1(\Sigma_g) / <c_1,\dots,c_n>\]
This will in turn give us a presentation for $\pi_1(X)$ by the standard generators $a_1,b_1,\dots,a_g,b_g$ with the relation $[a_1,b_1] \dots [a_g,b_g]$ along with relations given by each $c_i$ in terms of $a_1,b_1, \dots, a_g,b_g$.
An essential example is Matsumoto's fibration which is a genus $2$ Lefschetz fibration given by the positive factorization $(T_B T_C T_D T_\sigma)^2$ as given by the curves in figure ~\ref{fig:matsumoto}.
We then get a presentation for the fundamental group of the corresponding symplectic $4$-manifold on the generators $a_1,b_1,a_2,b_2$ with the relations
\[a_1a_2 = 1\]
\[a_1 \bar{b}_1a_2\bar{b}_2 = 1\]
\[b_2 b_1 = 1\]
\[[a_1,b_1] = 1\]
\[[a_1,b_1][a_2,b_2]=1\]
The first four relations come from the curves $B,C,D,\sigma$ and in that order.
We can now see that the fundamental group will be generated by $a_1$ and $b_1$ from the 1st and the 3rd relations.
Then from the fourth, we see that these two generators will commute.
All other relations will trivialize this, so we see that the fundamental group will isomorphic to the free abelian group on two generators.
One important observation will be that the presence of the curves $B,D,\sigma$ in our positive factorizations will imply that the fundamental group will be abelian and has presentation generated by $a_1, b_1$ with $a_1 = -a_2$, $b_1 = -b_2$.

Besides the fundamental group, we can also get other topological invariants from the positive factorization.
We will now bring our focus on to genus $2$ Lefschetz fibrations.
If a genus $2$ Lefschetz fibration has type $(n,s)$, then we can calculate many topological invariants from the type
\begin{prop}
Suppose a genus $2$ Lefschetz fibration has type $(n,s)$, then we can calculate the Euler characteristic $e$, signature $\sigma$, the holomorphic Euler characteristic $\chi_h$, and the Chern number $c_1^2$ by the formulas
\begin{itemize}
	\item $e = n+s - 4$
	\item $\sigma = -\frac{1}{5}(3n+s)$
	\item $\chi_h = \frac{1}{10}(n+2s)-1$
	\item $c_1^2 = \frac{1}{5}(n+7s) - 8$
\end{itemize}
\end{prop}
The Euler characteristic follows from the more general formula for a genus $g$ Lefschetz with $\ell$ critical points which is given by $e = 4-4g + \ell$.
For the signature, we take advantage of the fact that $\Gamma_2$ is hyperelliptic, i.e. every element of $\Gamma_2$ commutes with a fixed hyperelliptic involution.
Because of this, every genus $2$ positive factorization is fixed by a hyperelliptic involution and so we can apply Endo's signature formula to any genus $2$ Lefschetz fibration \cite{endo2000meyer}.
Then $\chi_h$ and $c_1^2$ can be calculated from the Euler characteristic and signature.
\subsection{Surgery on Lefschetz fibrations}
If we have Lefschetz fibrations given by positive factorizations, then certain manipulations of the positive factorizations will correspond to surgery on the Lefschetz fibrations.
Suppose we have two genus $g > 0$ Lefschetz fibrations $(X,f)$ and $(X',f')$ with positive factorizations $W$ and $W'$ respectively.
Then we can take a symplectic fiber sum of these fibrations along a non-singular fiber by some mapping class $\phi$.
This will preserve the fibration structure so we will still have a Lefschetz fibration which will have positive factorization $W^\phi W'$.
By Usher's work on symplectic fiber sums, this will be minimal \cite{usher2006minimality} (also see \cite{baykur2016minimality} where it was proven for fiber sums of Lefschetz fibrations).
\begin{prop}
Suppose we have two genus $g$ Lefschetz fibrations given by positive factorizations $W$ and $W'$ and $\phi \in \Gamma_g$.
Then the positive factorization $W^\phi W'$ will describe a minimal genus $g$ Lefschetz fibration.
\end{prop}
\noindent So we now have a method to produce many minimal Lefschetz fibrations by appending one positive factorization to another.

Rational blow downs were introduced by Fintushel and Stern and were shown by Symington to be possible in the symplectic category as well \cite{fintushel1995rational,symington2001generalized}.
In particular, we can blow down $(-4)$-spheres in a symplectic $4$-manifold.
If we have a Lefschetz fibration, this can be done via a monodromy substitution.
The lantern relation is a relation on the sphere with $4$ boundary components.
This is a relation between the product of the $4$ boundary twists and Dehn twists along three separating curves.
This can be generalized to $\Gamma_g$ by embedding the sphere with $4$ boundary components in $\Sigma_g$.
Suppose $a$ and $b$ are two simple closed curves in $\Sigma_g$ with geometric intersection number $2$ and algebraic intersection number $0$.
Then locally orient $a$ and $b$ at each intersection point so that each intersection point is positively oriented.
Resolve each intersection according to the local orientations according to the local orientation to get another simple closed curve $c$.
Then there is a neighborhood of $a \cup b$ will be a sphere with $4$ boundary components with four curves $d_1,d_2,d_3,d_4$ parallel to the boundary components and we have the relation
\[T_{d_1} T_{d_2} T_{d_3} T_{d_4} = T_a T_b T_c\]
If we have $4$ Dehn twists that bound a sphere in a positive factorization
\[W = \cdots T_{d_1}T_{d_2}T_{d_3}T_{d_4} \cdots\]
Then we can apply a lantern relation to get a new positive factorization
\[W' = \cdots T_a T_b T_c \cdots\]
We call this a lantern substitution and will be a crucial method for us to construct new positive factorizations because they correspond to a rational blow down.
\begin{prop} (Endo-Gurtas)

Suppose we have two positive factorizations $W$ and $W'$ describing genus $g$ Lefschetz fibrations $(X,f)$ and $(X',f')$.
If $W'$  is obtained from a lantern substitution from $W$, then $X'$ is the rational blow down of a symplectic $(-4)$-sphere in $X$ \cite{endo2010lantern}.
\end{prop}
\noindent This was later generalized to the rational blow downs of configurations $C_p$ of spheres via relations in the planar mapping class group \cite{endo2011monodromy}.
We will only need to rationally blow down a configuration $C_2$ in a genus $2$ Lefschetz fibration.

For a genus $2$ Lefschetz fibration, we can perform a lantern substitution when we have a subword $T_a^2 T_b^2$ where $a$ and $b$ are disjoint non-separating curves.
Then a lantern substitution will replace these $4$ non-separating Dehn twists with two non-separating Dehn twists and one separating Dehn twists.
Then a lantern substitution in a positive factorization of type $(n,s)$ will result in a positive factorization of type $(n-2,s+1)$.
\subsection{Small Lefschetz Fibrations and Exotic $4$-manifolds}
\begin{figure}
\centering
\includegraphics[scale = 0.7]{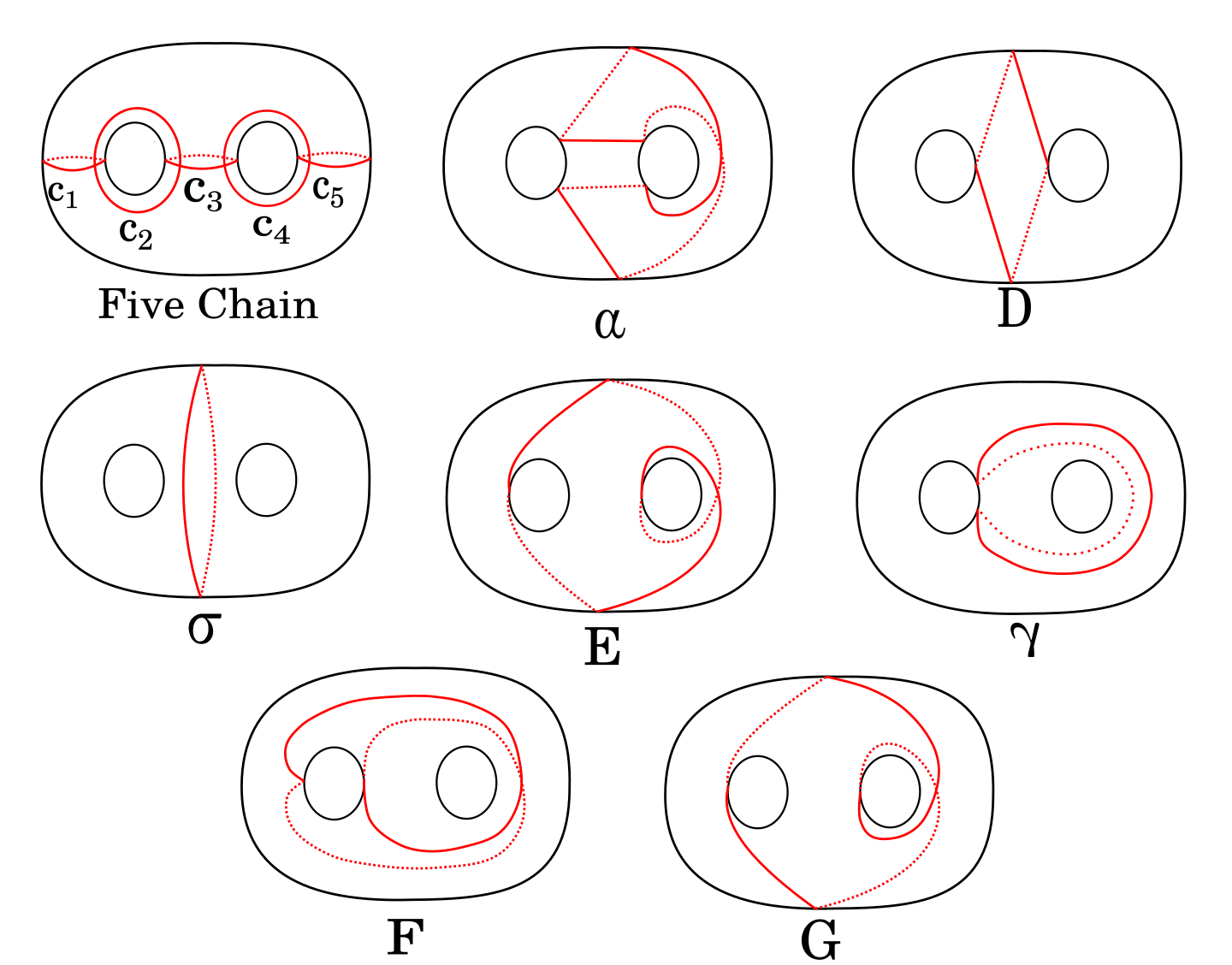}
\caption[Curve]{The five chain $c_1,c_2,c_3,c_4,c_5$, the generators $a_1,b_1,a_2,b_2$ for $\pi_1 (\Sigma_2)$, and the Dehn twist curves $B,C,D,\sigma,E,F,G, \alpha,\gamma$}
\label{fig:4_3_curves}
\end{figure}
Genus $0$ and $1$ Lefschetz fibrations are well understood, their total spaces are diffeomorphic to $\mathbb{CP}^2$, $\mathbb{CP}^1 \times \mathbb{CP}^1$, the elliptic surfaces $E(n)$, or blow ups of one of these spaces.
One is then led to ask about genus $2$ Lefschetz fibrations.
Baykur and Korkmaz addressed this question for simply connected, minimal genus $2$ Lefschetz fibrations with the additional restriction that the underlying symplectic $4$-manifold is small, i.e. has $b^+ \le 3$.
They were able to establish the following
\begin{thm} (Baykur-Korkmaz)

Any simply connected, minimal genus $2$ Lefschetz $(X,f)$ with $b^+ \le 3$ is homeomorphic to $\mathbb{CP}^2 \# p \overline{\mathbb{CP}}^2$ for some $7 \le p \le 9$ or $3\mathbb{CP}^2 \# q \overline{\mathbb{CP}}^2$ for some $11 \le q \le 19$.
Furthermore, there are positive factorizations of types $(8,6),(10,5),(12,4)$ and $(12,9),(14,8),(16,7),(18,6),(20,5),(22,4),(24,3),(26,2)$ that realize minimal, decomposable genus $2$ Lefschetz fibrations whose total spaces are $\mathbb{CP}^2 \# p \overline{\mathbb{CP}}^2$ for $p = 7,8,9$ and $3\mathbb{CP}^2 \# q \overline{\mathbb{CP}}^2$ for $q = 12, \dots, 19$ respectively. \cite{baykur2017small}
\end{thm}
The main tool to establish this was the use of fiber sums of small positive factorization.
To achieve this they first needed to construct a positive factorization $W$ of type $(4,3)$ which describes the smallest genus $2$ Lefschetz fibration (small as in fewest number of critical points).
This was found by taking the positive factorization $(T_{c_1}T_{c_2}T_{c_3}T_{c_4})^{10}$ of type $(40,0)$ and finding three $2$-chain substitutions.
By taking fiber sums of $W$ along with other small positive factorizations, they were able to construct all of the spaces in theorem $2$.
We will not use the positive factorization as given in \cite{baykur2017small}, but instead use an alternate $(4,3)$ factorization communicated to the author by Noriyuki Hamada and will appear in \cite{baykurhamada18}.
\begin{thm} 
The positive factorization given by the curves in figure ~\ref{fig:4_3_curves}
\[W = T_\alpha T_D T_\sigma T_E T_\gamma T_F T_G\]
prescribes the smallest genus $2$ Lefschetz fibration \cite{baykurhamada18,baykur2017small}.
\end{thm}
This alternate factorization has much simpler curves and for this reason we will use this alternate factorization for our constructions.
\section{Positive Factorization of Type $(10,10)$}
We will use three copies of the alternate $(4,3)$ factorization and a lantern substitution to construct our $(10,10)$ positive factorization.
Here we will write the alternate $(4,3)$ factorization as $W = \alpha D \sigma E \gamma FG$ and in general will denote a Dehn twist along a curve $c$ instead of writing $T_c$.
The inverse of the Dehn twist $c$ will be denoted by $\bar{c}$.
Also, we will denote Dehn twists along non-separating curves by capital letters and Dehn twist along separating curves by greek letters.
Furthermore, if $\phi$ is some mapping class, then we let ${}_\phi c$ denote the Dehn twist along $\phi (c)$ which is equal to $\phi c \bar{\phi}$.
Also, the standard $5$-chain $c_1,c_2,c_3,c_4,c_5$ will simply be denoted by the numbers $1,2,3,4,5$.

First we will take a twisted fiber sum of two copies of the alternate $(4,3)$ positive factorization and apply Hurwitz moves to get a $(6,8)$ positive factorization.
Take a twisting of $W$ by the Dehn twist $4$ which we write as
\[W^4 = {}_4\alpha {}_4D \sigma {}_4E {}_4 \gamma F {}_4G\]
We can then move ${}_4\alpha$ to the end of this positive factorization since this composition of Dehn twists is equal to the identity we can conjugate by ${}_4\bar{\alpha}$ and it will still be the identity (one can also use Hurwitz moves to move this if one wants to preserve the fibration which will not be necessary for us).
\[W' = {}_{4} D \sigma {}_{4} E \gamma F {}_{4} G {}_{4} \alpha\]
By the same logic we can bring $\alpha$ to the end of $W$ to get $D \sigma E \gamma F G \alpha$ and we can then insert $W^4$ into this factorization like so
\[D W^4 \sigma E \gamma F G \alpha = D {}_{4} D \sigma {}_{4} E \gamma F {}_{4} G {}_{4} \alpha \sigma E \gamma F G \alpha\]
Then we can apply a Hurwitz move to push $D$ past ${}_{4} D$
\[{}_{D4} D D \sigma {}_{4} E \gamma F {}_{4} G {}_{4} \alpha \sigma E \gamma F G \alpha\]
Since $D$ and $4$ intersect geometrically once, we can apply the braid relation to ${}_{D4} D=D4D\bar{4}\bar{D} = 4D4\bar{4}\bar{D}=4$.
So we have:
\[4 D \sigma {}_{4} E \gamma F {}_{4} G {}_{4} \alpha \sigma E \gamma F G \alpha\]
By Hurwitz moving the two $F$ Dehn twists to the end, this is Hurwitz equivalent to:
\[4 D \sigma {}_{4} E \gamma {}_{F4} G {}_{F4} \alpha {}_{F} \sigma {}_{F} E {}_{F} \gamma {}_{FF} G {}_{FF} \alpha FF\]
We can then move the two $F$ Dehn twists to the front to get
\[V:= FF 4 D \sigma {}_{4} E \gamma {}_{F4} G {}_{F4} \alpha {}_{F} \sigma {}_{F} E {}_{F} \gamma {}_{FF} G {}_{FF} \alpha\]
Now take a copy of $D \sigma E \gamma F G \alpha$ and twist it by $1 \bar{4}$ to get
\[U := {}_{4} D \sigma {}_{\bar{1} 4} E \gamma {}_{\bar{1}} F {}_{\bar{1} 4} G {}_4 \alpha\]
Take a fiber sum of $U$ and $V$ like so to get a positive factorization of type $(12,9)$ on some minimal symplectic $4$-manifold $Y'$
\[U (FF 4 D \sigma {}_{4} E \gamma {}_{F4} G {}_{F4} \alpha {}_{F} \sigma {}_{F} E {}_{F} \gamma {}_{FF} G {}_{FF} \alpha)\]
We can then Hurwitz move the Dehn twists in $U$ past $FF4D$ to get
\[({}_UF{}_UF {}_U4 {}_UD) U (\sigma {}_{4} E \gamma {}_{F4} G {}_{F4} \alpha {}_{F} 
\sigma {}_{F} E {}_{F} \gamma {}_{FF} G {}_{FF} \alpha\]
However, since $U$ is a positive factorization it will act by the identity on these Dehn twists.
We also expand $U$ out to get 
\[X' = FF 4 D ({}_{4} D \sigma {}_{\bar{1} 4} E \gamma {}_{\bar{1}} F {}_{\bar{1} 4} G {}_4 \alpha)\sigma {}_{4} E \gamma {}_{F4} G {}_{F4} \alpha {}_{F} \sigma {}_{F} E {}_{F} \gamma {}_{FF} G {}_{FF} \alpha\]
We want to take a moment to emphasize why we adjoined $U$ and $V$ first and then applying Hurwitz moves instead of just inserting $U$ into $V$.
When we adjoin two positive factorizations, the corresponding Lefschetz fibration will be a fiber sum and be minimal.
Then when we apply these Hurwitz moves, we get a positive factorization which has an isomorphic Lefschetz fibration which therefore must be minimal.
By simply inserting $U$ into $V$ it is not clear that the resulting positive factorization will be equivalent to a fiber sum.
After our digression to observe minimality, we are ready to find a lantern relation in our positive factorization.
Hurwitz move $D$ past ${}_4 D$ to get
\[FF 4  {}_{D4} D D \sigma {}_{\bar{1} 4} E \gamma {}_{\bar{1}} F {}_{\bar{1} 4} G {}_4 \alpha \sigma {}_{4} E \gamma {}_{F4} G {}_{F4} \alpha {}_{F} \sigma {}_{F} E {}_{F} \gamma {}_{FF} G {}_{FF} \alpha\]
\begin{figure}
\centering
\includegraphics[scale=1.0]{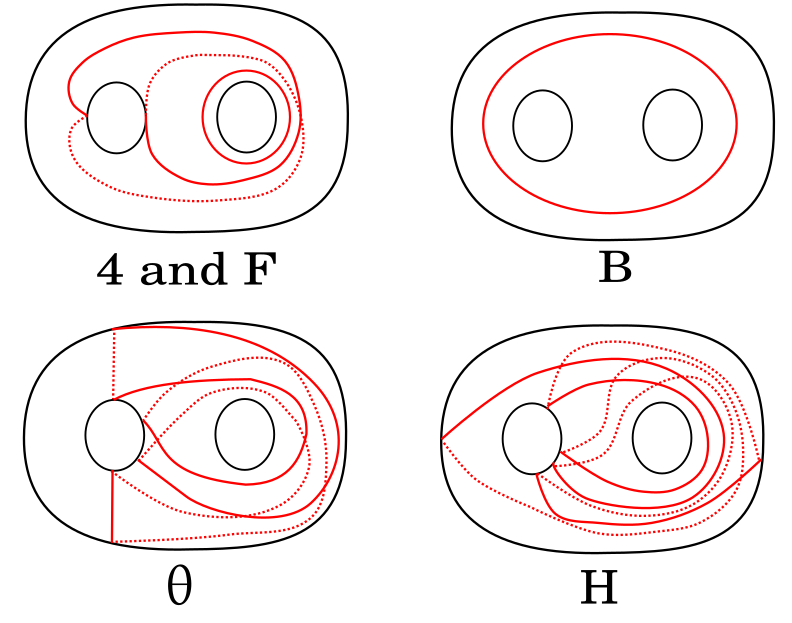}
\caption[Lantern]{Lantern Relation for $4$ and $F$}
\label{fig:lantern}
\end{figure}
As we observed earlier, ${}_{D4}D=4$ so we get the following positive factorization that contains the lantern $FF44$
\[FF44 D \sigma {}_{\bar{1} 4} E \gamma {}_{\bar{1}} F {}_{\bar{1} 4} G {}_4 \alpha \sigma {}_{4} E \gamma {}_{F4} G {}_{F4} \alpha {}_{F} \sigma {}_{F} E {}_{F} \gamma {}_{FF} G {}_{FF} \alpha\]
Then we can apply a lantern substitution to $FF44$ and get a positive factorization $X$ of type $(10,10)$.
The lantern relation here will be given by the curves in figure ~\ref{fig:lantern} as $44FF = B \theta H$
\[X := B \theta H D \sigma {}_{\bar{1} 4} E \gamma {}_{\bar{1}} F {}_{\bar{1} 4} G {}_4 \alpha \sigma {}_{4} E \gamma {}_{F4} G {}_{F4} \alpha {}_{F} \sigma {}_{F} E {}_{F} \gamma {}_{FF} G {}_{FF} \alpha\]
I claim that the corresponding Lefschetz fibration is on a simply connected, minimal symplectic $4$-manifold which we denote by $Y$.

To show simple connectedness, we have a presentation for the corresponding symplectic $4$-manifold as the quotient of $\pi_1 (\Sigma_2) \cong <a_1,a_2,b_1,b_2 |[a_1,b_1][a_2,b_2]>$ by the normal subgroup generated by curves in the positive factorization.
First observe as we did earlier for Matsumoto's fibration that the presence of the Dehn twists $B,D,\sigma$ imply that the fundamental group is abelian and generated by $a_1,b_1$ since $a_2 = -a_1$ and $b_2 = -b_1$.
Since the fundamental group is abelian, we only need to consider the homology classes of the remaining curves.
Now observe that ${}_{\bar{1}4}G$ is equal to ${}_5 B$ as curves, then by applying the Picard Lefschetz theorem we see that
\[[{}_{\bar{1}4}G] = [{}_5 B] = [5] + \hat{i}(5,B)[B] = a_1 + a_2 + b_2\]
Next if we calculate the algebraic intersection numbers of $F$ with the basis $a_1,b_1,a_2,b_2$ we can to get the homology class of $F$.
We have that $\hat{i}(a_1,F) = \hat{i}(b_2,F) = 0$, $\hat{i}(b_1,F) = 1$, and $\hat{i}(b_2,F) = 2$.
This enough to conclude that $[F] = a_1 + 2a_2$ and applying Picard-Lefschetz again, we see that
\[[{}_{\bar{1}}F] = [F] - \hat{i}(1,F)[1] = a_1 + 2 a_2 - b_1\]
Now we have enough information to conclude that the space corresponding to our $(10,10)$ positive factorization is simply connected.
As we observed for Matsumoto's fibration, the presence of the Dehn twists $B,D,\sigma$ imply that the fundamental group is abelian and that $a_1=-a_2$, $b_1 = -b_2$ in our presentation for $\pi_1(Y)$.
With the relation $a_1+a_2+b_2 = 0$ coming from ${}_{\bar{1}4}G$ and that $a_1 = -a_2$, we see that the relation $b_2 = 0$ holds.
Then since $b_1 = -b_2$, we have that $b_1 = 0$.
Then the presence of ${}_{\bar{1}}F$ implies that $a_1+2a_2-b_1=0$ and substituting $a_1=-a_2$ into this implies that $a_2 = b_1$.
Therefore $a_2 = 0$ since $b_1= 0$ and together with $a_1 =-a_2$, we can conclude $a_1 = 0$.
We have now seen that all generators are equal to zero in the fundamental group, therefore it is trivial.

We now move on to conclude minimality of this Lefschetz fibration.
We took a lantern substitution of a positive factorization that is Hurwitz equivalent to $UV$ which represents a minimal Lefschetz fibration.
Endo and Gurtas showed that the lantern substitutions in positive factorizations corresponds to a rational blow-down of a $(-4)$-sphere in the corresponding symplectic $4$-manifolds \cite{endo2010lantern}.
Thus the symplectic $4$-manifold coming from the positive factorization $X$ is a rational blow-down of a minimal symplectic $4$-manifold $Y'$ along a $(-4)$-sphere.
However, the rational blow-down of a $(-4)$-sphere will be minimal if there are no exceptional spheres in the complement of the $(-4)$-sphere and there is not a pair of disjoint exceptional spheres each intersecting the $(-4)$-sphere transversely at a single point \cite{dorfmeister2013minimality}.
However, since the Lefschetz fibration we are rationally blowing down is minimal, this condition must hold.
Therefore, the symplectic $4$-manifold $Y$ corresponding to the positive factorization $X$ is minimal.

Using the formulas we listed earlier and simple connectedness, we can conclude that $Y$ has $b_2 = 14$ and $\sigma = -8$.
Since the signature is not divisible by $16$, this manifold has odd intersection form and therefore by Freedman's classification theorem we can conclude that $Y$ is homeomorphic to $3 \mathbb{CP}^2 \# 11 \overline{\mathbb{CP}}^2$ \cite{freedman1982topology}.
However, $Y$ won't be diffeomorphic to $3 \mathbb{CP}^2 \# 11 \overline{\mathbb{CP}}^2$ since the latter is minimal while the former is not.
\begin{thm}
There exist a genus $2$ Lefschetz fibration whose total space is a minimal symplectic $4$-manifold which is homeomorphic, but not diffeomorphic to $3 \mathbb{CP}^2 \# 11 \overline{\mathbb{CP}}^2$
\end{thm}
In particular, this completes Baykur and Korkmaz's geography of simply connected minimal small genus $2$ Lefschetz fibrations.
\begin{remark}
The main idea of this construction is that the positive factorization $V$ is a ``fake $(6,7)$".
Namely it contains the subword $4FF$ which consists of three of the curves needed for a lantern substitution.
If $T$ is some positive factorization of type $(n,s)$, then by choosing an appropriate twisting $\phi$, one can guarantee that the fiber sum $T^\phi V$ has the fourth curve needed for a lantern substitution.
Then $T^\phi V$ has type $(n+8,s+6)$ and after applying a lantern substitution we get a positive factorization of type $(n+6,s+7)$.
This idea will be exploited in the next section to construct a family of interesting Lefschetz fibrations. 
\end{remark}
\section{Geography Findings}
\subsection{Restriction on the Geography of Genus $2$ Lefschetz fibrations}
The geography problem for a given class of $4$-manifolds is to determine what pairs of holomorphic Euler characteristic and first Chern numbers can be realized by this class of manifolds.
Recall that if we have a genus $2$ Lefschetz fibration of type $(n,s)$ we have the following formulas for the holomorphic Euler characteristic signature and the first Chern number.
\[c_1^2 = \frac{1}{5}(n+7s)-8 = 2m + s - 8\]
\[\chi_h = \frac{1}{10}(n+2s)-1 = m-1\]
Where $m = \frac{1}{10}(n+2s)$ which is always an integer.
It then follows that $c_1^2 = 2\chi_h + s - 6 \ge 2\chi_h - 6$ by substituting $m = \chi_h + 1$ and using that $s \ge 0$.
We have now shown that the Noether inequality holds for genus $2$ Lefschetz fibrations.

The next restriction we will observe on $(c_1^2,\chi_h)$ will follow from Baykur and Korkmaz's observation that $b^- \ge s+1$ for a Lefschetz fibration of type $(n,s)$.
Combining this inequality and Poincare duality that the Betti numbers satisfy $b_0=b_4=1$ and $b_1=b_3$
\[e - \sigma = 2-2b_1 + 2b^- \ge 2-2b_1 + 2(s+1)\]
Substituting our formulas for the Euler characteristic and signature in terms of the type $(n,s)$
\begin{align*}
(n+s-4) + \frac{1}{5}(3n+s) \ge 2-2b_1 + 2(s+1) &\Rightarrow 8(\frac{n+2s}{10}) \ge 4 -b_1 + 2s \\
							&\Rightarrow 8m \ge 4-b_1 + 2s 
\end{align*}
We now use our formulas for $\chi_h$ and $c_1^2$ in terms of $m$ and $s$ to derive the following inequality
\begin{align*}
8m \ge 4-b_1 + 2s 	&\Rightarrow 12m -16 \ge 4 - b_1 + 2s + 4m -16 \\
					&\Rightarrow 12(\chi_h+1) -16 \ge 4 - b_1 + 2c_1^2 \\
					&\Rightarrow c_1^2 \le 6\chi_h - 4 + \frac{b_1}{2}
\end{align*}
Since the fundamental group of a genus $2$ Lefschetz fibration will be isomorphic to a  quotient of $\pi_1 (\Sigma_2)$ by a subgroup containing at least one non-separating curve (any non-trivial genus $2$ Lefschetz fibration has $n > 0$), then $b_1 \le 3$ and we can conclude that $c_1^2 \le 6\chi_h - \frac{5}{2}$.
Since $c_1^2$ and $\chi_h$ are both integers we can conclude that $c_1^2 \le 6\chi_h - 3$.
For a simply connected genus $2$ Lefschetz fibration, we have $b_1=0$ and therefore $c_1^2 \le 6\chi_h - 4$ which is equivalent to the strict inequality $c_1^2 < 6\chi_h - 3$ since $c_1^2$ and $\chi_h$ are both integers.
\begin{thm}
Genus $2$ Lefschetz fibrations satisfy the Noether inequality $c_1^2 \ge 2\chi_h - 6$.
Furthermore, they satisfy the inequality $c_1^2 \le 6\chi_h - 3$ and this inequality can be made strict for simply connected genus $2$ Lefschetz fibrations.
\end{thm}
\noindent We now would like to construct positive factorizations that realize genus $2$ Lefschetz fibrations.
In particular, we will focus on constructing examples that are simply connected and minimal.
\begin{remark}
One can't help but draw an analogy between this result and the geography of complex surfaces of general type.
Both satisfy the Noether inequality and we think of the inequality $c_1^2 \le 6\chi_h - 3$ for genus $2$ Lefschetz fibrations as the analogues of the BMY inequality for complex surfaces of general type.
This analogy will become stronger with and will guide the following constructions.
\end{remark}
\begin{remark}
One is immediately tempted to ask if a similar result holds for higher genus Lefschetz fibrations.
This seems plausible, but difficult due to the loss of a restriction like $n+2s = 0 \ (mod \ 10)$ on the type like we had for genus $2$ Lefschetz fibrations since the abelianization of $\Gamma_g$ is trivial for $g > 2$.
Furthermore, we lose hyperellipticity of $\Gamma_g$ for $g > 2$ and so there will be positive factorizations that are not hyperelliptic.
This means we will not be able to apply Endo's signature formula to all Lefschetz fibrations with genus greater than $2$.
Without Endo's signature formula, it becomes very difficult to deduce any restrictions on the topological invariants of higher genus Lefschetz fibrations.
\end{remark}
\subsection{Below the BK-line}
First we will show that the techniques of Baykur and Korkmaz will realize many values of $(\chi_h,c_1^2)$.
Taking fiber sums of $(4,3)$ positive factorizations will allow one to construct positive factorizations of type $(4k,3k)$ that describe simply connected minimal genus $2$ Lefschetz fibrations for $k \ge 2$.
Using the formulas we had earlier, we can calculate that these will have $\chi_h = k-1$ and $c_1^2 = 5k-8$.
These satisfy the linear relation $c_1^2 = 5\chi_h-3$ which we will refer to as the BK-line since fiber sums of Baykur-Korkmaz $(4,3)$ positive factorization will realize all $(\chi_h,c_1^2)$ on this line and it also represents the upper limit of Baykur and Korkmaz's technique of taking fiber sums to find simply connected genus $2$ Lefschetz fibrations.
We restate their results in the following form
\begin{thm}(Baykur-Korkmaz)

For $(x,y)$ with $x=1,2$ and $0 \le y \le 5x-3$, there exists minimal, simply connected genus $2$ Lefschetz fibrations such that $\chi_h = x$ and $c_1^2 = y$.
\end{thm}
We will observe that this results can be extended for $(x,y)$ with $2x-6 \le y \le 5x-3$, $y \ge 0$ using fiber sums.
The strategy to accomplish this is to fix a number of separating cycles $s$ and construct positive factorizations describing simply connected, minimal genus $2$ Lefschetz fibrations with this number of separating cycles.
One can observe from our formulas for the first Chern number $c_1^2$ and $\chi_h$ for a genus $2$ Lefschetz fibration of type $(n,s)$ satisfy
\[c_1^2 = \frac{1}{5}(n+7s)-8 = 2(\frac{n+2s}{10} - 1) -6 + s = 2\chi_h -6 + s\]
Therefore, a genus $2$ Lefschetz fibration with $s$ separating cycles lays on the line $c_1^2 = 2\chi_h - 6 + s$.
We observe that when $s = 0$ we recover the Noether line as our lower bound for the geography of genus $2$ Lefschetz fibration.
For this reason, we call this line the $s$-Noether line.

We will only need to construct the smallest two positive factorizations on the $s$-Noether line that are below the BK-line with positive $c_1^2$ that describe simply connected, minimal genus $2$ Lefschetz fibrations.
These two will be given by positive factorization of type $(n,s)$ and $(n+10,s)$.
Then by taking fiber sums of these with $(20,0)$ positive factorization $(1234554321)^2$ to get positive factorizations of $(n+10k,s)$.
For example, when $s = 0$ the smallest positive factorization which will have positive $c_1^2$ will have type $(40,0)$.
A positive factorization with this type will be $(1234554321)^4$ which gives a simply connected minimal genus $2$ Lefschetz fibration.
This will be minimal since it is a fiber sum of two copies of the positive factorization $(1234554321)^2$ and simply connected since the curves $c_1,c_2,c_3,c_4,c_5$ will kill the fundamental group.
The positive factorization $(1234554321)^2 (12345)^6$ will be of type $(50,0)$ and will also be simply connected and minimal.
Then by taking fiber sums of these two positive factorizations with copies of $(1234554321)^2$ will realize minimal simply connected genus $2$ Lefschetz on all point of the Noether line with $c_1^2 \ge 0$.

For the $s=1$ case, the smallest two positive factorizations would be of types $(38,1)$ and $(48,1)$.
These can be realized by fiber sums of the positive factorization $W_4$ of type $(18,1)$ described by Baykur and Korkmaz with the positive factorizations $(1234554321)^2$ and $(12345)^6$.
For the $s=2$, the smallest two positive factorizations will be of types $(26,2)$ and $(36,2)$.
These can be realized by Matsumoto's $(6,2)$ positive factorization fiber summed with $(1234554321)^2$ and $(12345)^6$.
For the $s = 3$ case, the smallest positive positive factorizations will be of types $(24,3)$ and $(34,3)$.
These can be realized by fiber sums of Baykur and Korkmaz's $(4,3)$ positive factorization with $(1234554321)^2$ and $(12345)^6$.
Then by taking repeated fiber sums of these positive factorizations with $(1234554321)^2$ one can get the entirety of the $s$-line for $s=1,2,3$.

If a positive factorization of type $(n,s)$ is under the $BK$-line, then we can substitute our formulas for $c_1^2$ and $\chi_h$ in terms of $(n,s)$ into the inequality $c_1^2 \le 5\chi_h - 3$.
These will simplify to $4s \le 3n$ and let $n(s)$ be the smallest positive integer such that this inequality is satisfied for a given $s$ and $n+2s = 0 \ (mod \ 10)$.
Then observe that $n(s+3) = n(s)+4$ and therefore if we have positive factorizations of type $(n(s),s)$ and $(n(s)+10,s)$ then we can simply take a fiber sum with the Baykur-Korkmaz $(4,3)$ positive factorization to get $(n(s+3),s+3)$ and $(n(s+3)+10,s+3)$.
So we can reduce this to the $s=4,5,6$ cases.
Fortunately, positive factorizations that describe simply connected, minimal genus $2$ Lefschetz fibrations were given by Baykur and Korkmaz.
For $s = 4$, these are their $\widetilde{W}_3$ and $\widehat{W}_6$; for $s = 5$, these are their $\widetilde{W}_2$ and $\widehat{W}_6$; for $s=6$, these are their $\tilde{W}_1$ and $\widehat{W}_4$.
Now by taking repeated fiber sums with the $(4,3)$ positive factorizations and $(1234554321)^2$ we can realize all $(\chi_h,c_1^2)$ below the BK-line.
\begin{thm}
For $(x,y)$ such that $2x-6 \le y \le 5x-3$ and $y \ge 0$, there exists minimal, simply connected genus $2$ Lefschetz fibrations such that $\chi_h = x$ and $c_1^2 = y$.
\end{thm}
\subsection{Above the BK-line}
The positive factorization $X$ of type $(10,10)$ described earlier has $(\chi_h,c_1^2) = (2,8)$ which is above the $BK$-line.
By taking fiber sums of $X$ with $(4,3)$ positive factorizations one can get positive factorizations of type $(10+4k,10+3k)$ which will have $(\chi_h,c_1^2) = (2+k,5k+8)$.
These will be all be minimal and simply connected and will have $(\chi_h,c_1^2)$ with $c_1^2 = 5\chi_h-2$.
This will realize an infinite number of distinct genus $2$ Lefschetz fibrations, but they stay very close to the BK-line.
To realize genus $2$ Lefschetz fibrations that are farther away from the $BK$-line, we will construct a family of positive factorizations $X(t)$ with $X(1) = X$ and $X(t)$ has type $(4+6t,3+7t)$.
These will have $(\chi_h,c_1^2) = (2t,11t-3)$ which satisfy the relation $c_1^2 = 5.5\chi_h -3$.
Then we can fill the region between $5\chi_h < c_1^2 \le 5.5\chi_h - 3$ by partitioning the region by the lines $c_1^2 = 5\chi_h - 3 + t$ for $t$ a non-negative integer.
We call this the $BK(t)$-line since for $t = 0$ we recover the $BK$-line and for $t > o$ we get lines that run parallel to the $BK$-line.

First we discuss how to construct the $X(t)$ positive factorizations inductively by fiber sums and lantern substitutions.
Let $X(1) = X$ and observe that the presence of the Dehn twists $B,D,\sigma, {}_{\bar{1}}F, {}_{\bar{1}{4}}G$ trivialize the fundamental group.
Suppose we have constructed $X(t)$ such that we have these Dehn twists present and they appear in this order.
Furthermore, suppose that $X(t)$ has the desired type of $(4+6t,3+7t)$
We can then Hurwitz move these to the front of $X(t)$ and observe that there will be a non-separating Dehn twist that we can Hurwitz move to the back.
So we can use Hurwitz move to change $X(t)$ of the form
\[BD\sigma {}_{\bar{1}}F {}_{\bar{1}{4}}G \cdots J\]
Then there is some mapping class $\phi$ that sends $4$ to the Dehn twist $J$ and we can fiber sum $V^\phi$ with the above to get
\[(BD\sigma {}_{\bar{1}}F {}_{\bar{1}{4}}G \cdots J) V^\phi= BD {}_{\bar{1}}F {}_{\bar{1}{4}}G, \sigma \cdots J {}_\phi F{}_\phi F {}_\phi 4 \cdots\]
Then ${}_\phi 4 = J$ and therefore we have
\[J {}_\phi F {}_\phi F {}_\phi 4 = {}_\phi 4{}_\phi F {}_\phi F {}_\phi 4 = {}_\phi (4FF4)\]
Since $4$ and $F$ are disjoint, twisting them by $\phi$ they will be disjoint.
Therefore we can apply a lantern substitution to get 
\[X(t+1) := BD\sigma {}_{\bar{1}}F {}_{\bar{1}{4}}G \cdots {}_\phi B {}_\phi \theta {}_\phi H \cdots\]
This positive factorization will describe a simply connected genus $2$ Lefschetz fibration due to the presence of the Dehn twists $B,D, \sigma {}_{\bar{1}}F, {}_{\bar{1}{4}}G$.
Then since we did a lantern substitution, $X(t+1)$ will describe a genus $2$ Lefschetz fibration that is a rational blow down of a $(-4)$-sphere.
However, the rational blow down will be on the fiber sum of the Lefschetz fibration given by the positive factorizations $V$ and $X(t)$.
This will be a minimal genus $2$ Lefschetz fibration since it is a fiber sum.
We can now conclude that $X(t+1)$ is minimal since it is the blow down of a $(-4)$-sphere \cite{dorfmeister2013minimality}.
If $X(t)$ has type $(4+6t,3+7t)$, then $X(t+1)$ will have the desired type of $(4+6(t+1),3+7(t+1))$ since $V$ has type $(8,6)$ and the lantern substitution reduce the non-separating cycles by $2$ and increases the number of separating cycles by $1$.
We have now constructed the desired $X(t)$ of the desired type which will be simply connected and minimal.

The $X(t)$ will then give Lefschetz fibrations with $c_1^2 = 5.5\chi_h - 3$.
To fill the region $5\chi_h - 3 < c_1^2 \le 5.5\chi_h - 3$ we use the partition into $BK(t)$ lines and give positive factorizations that describe Lefschetz fibrations on this line.
The $BK(t)$ line will intersect the line $c_1^2 = 5.5\chi_h - 3$ at the point $(2t,11t-3)$.
Then $X(t)$ will realize a simply connected, minimal genus $2$ Lefschetz fibration at this point.
We can take fiber sums of $X(t)$ with $k$ copies of the $(4,3)$ positive factorization.
This will have type $(4 + 6t + 4k,3 + 7t + 3k)$ and we can calculate that $(\chi_h,c_1^2) = (2t+k,11t+5k-3)$.
These will realize all points on the $BK(t)$-line which will be minimal and simply connected.
The positive factorization will still have the Dehn twists $B,D,\sigma, {}_{\bar{1}}F, {}_{\bar{1}{4}}G$.
They will be fiber sums and will therefore by minimal.
We can now conclude the following
\begin{thm}
For $(x,y)$ such that $5x-3 < y \le 5.5x-3$ and $y \ge 0$, there exists minimal, simply connected genus $2$ Lefschetz fibrations such that $\chi_h = x$ and $c_1^2 = y$.
\end{thm}
\begin{remark}
The positive factorization $X(2)$ has type $(16,17)$ and as far as the author is aware of, is the smallest positive factorization with more separating Dehn twists than non-separating that has been constructed in the literature.
The smallest types with this property that are not ruled out by theorem $5$ for example would be $(6,7)$ and $(8,11)$.
The only smaller type that has this property and does not have $c_1^2 = 6\chi_h-3$ would be the type $(12,14)$.
\end{remark}
Combining this with theorem $7$ and $8$, we can realize all $(x,y)$ with $2x-6 \le y \le 5.5x-3$ and $y \ge 0$.
Seeing how we construct positive factorizations in steps that go from the bound $c_1^2 = 5\chi_h -3$ to the bound $c_1^2 - 5.5\chi_h - 3$, and that we have the bound $c_1^2 \le 6\chi_h - 3$ we are interested in constructing positive factorizations such that $c_1^2 = a\chi_h -3$ with $a$ close to $6$.
Such $a$ we call the slope of the positive factorization or genus $2$ Lefschetz fibration.
This will be bounded above by $6$ and will be equal to $6$ with positive factorizations of type $(n,s)$ with $s=2n - 5$.

\begin{remark}
With the results established so far on the geography of genus $2$ Lefschetz fibrations, the only Lefschetz fibration on the line $c_1^2 = 6\chi_h - 3$ given by an explicit positive factorization are Baykur-Korkmaz's $(4,3)$ positive factorization.
Xiao has constructed Lefschetz fibrations of type $(6,7)$ and $(12,19)$ which would also lay on this line \cite{xiao2006surfaces}.
However, positive factorizations for these Lefschetz fibrations are not known.
Even Xiao's Lefschetz fibration of type $(4,3)$ is not known to agree with the Lefschetz fibration given by Baykur and Korkmaz's $(4,3)$ positive factorization despite being on the same underlying manifold.
\end{remark}

Due to the lack of positive factorizations which describe Lefschetz fibrations with $c_1^2 = 6\chi_h - 3$, one might instead ask how close can we get to this line using positive factorizations.
If a genus $2$ Lefschetz fibration has $c_1^2 = a\chi_h - 3$, we call $a$ the slope of the Lefschetz fibration.
The slope will only depend on the type $(n,s)$
\[a = \frac{c_1^2 + 3}{\chi_h} = 2 + \frac{10s-30}{n+2s-10}\]
This will equal $6$ for genus $2$ Lefschetz fibrations with $c_1^2 = 6\chi_h - 3$.
Theorems $4$ and $5$ will realize all rational slopes in the interval $[2,5.5]$.
Due to the lack of positive factorizations describing genus $2$ Lefschetz fibrations with $c_1^2 = 6\chi_h-3$ or equivalently having type $(n,2n-5)$, we instead try to find positive factorizations with large slope.
This will result in Lefschetz fibrations close to the line $c_1^2 = 6\chi_h - 3$.

The largest slope we can realize is given by a positive factorization of type $(24,38)$ and will have slope $\frac{53}{9} \approx 5.9$.
To construct this positive factorization, we find two $2$-chain relations in a positive factorization of type $(48,36)$.
We can insert a copy of the $(4,3)$ positive factorization $D \sigma E \gamma FG \alpha$ into $U = {}_{4} D \sigma {}_{\bar{1} 4} E \gamma {}_{\bar{1}} F {}_{\bar{1} 4} G {}_4 \alpha$ like so:
\[{}_{4} D (D \sigma E \gamma FG \alpha)\sigma {}_{\bar{1} 4} E \gamma {}_{\bar{1}} F {}_{\bar{1} 4} G {}_4 \alpha\]
We can then Hurwitz move $F$ and ${}_{\bar{1}}F$ to the end of this positive factorization
\[{}_4 D D Z' F {}_{\bar{1}}F\]
Where $Z'$ is some subword of type $(4,6)$.
We can then take six copies of this to get
\[({}_4 D D)^6(Z' F {}_{\bar{1}}F)^6\]
We can then Hurwitz move the six copies of ${}_{\bar{1}}F F$ to the end of this positive factorization to get
\[({}_4 D D)^6 Z (F {}_{\bar{1}}F)^6\]
Where $Z$ is a subword of type $(24,36)$.
We can then apply the $2$-chain relation to $({}_4 D D)^6$ and $(F {}_{\bar{1}}F)^6$ and replace these with some separating curves $\delta$ and $\omega$
\[\delta Z \omega\]
This will be a positive factorization of type $(24,38)$.
It is not necessarily simply connected or minimal, however it is of interest due to its large slope
\begin{thm}
There exists a genus $2$ Lefschetz fibration of type $(24,38)$ given by an explicit positive factorization.
\end{thm}
\begin{remark}
One can play a similar game to the one we did in theorem $8$ and get positive factorizations $M(t) = W (\delta Z \omega)^t$ which will have type $(4 + 24t,3 + 38t)$.
This will have $(\chi_h,c_1^2) = (10t,58t-3)$ and each of these will then be on the line $c_1^2 = 5.8\chi_h - 3$ and have slope $5.8$.
Then one can take fiber sums of these $M(t)$ with $k$ copies of the $(4,3)$ positive factorization to get positive factorizations with type $(4+24t+4k,3+38t+3k)$.
One can then show that these will have slope $5.8 - .8\frac{k}{10t+k}$ and by varying these $k$ and $t$, one can realize a set of slopes that will be dense in the interval $[5.5,5.8]$.
It would be interesting to get a similar result to Roulleau and Urz{\'u}a construction that realizes Chern slopes densely in the interval $[2,3]$ by finding genus $2$ Lefschetz fibrations with slopes dense in the interval $[5,6]$ or arbitrarily close to $6$ \cite{roulleau2015chern}.
\end{remark}
\section*{Acknowledgments}
I would like to thank my undergraduate thesis advisor R. Inanc Baykur for his guidance through this project and my last year at UMass Amherst.
Without his help and patience, this would not have been possible.

%


Department of Mathematics, University of Texas, Austin, Texas 78712-1202, USA

\textit{Email Address}: kainakamura@utexas.edu
\end{document}